\documentclass[a4paper,10pt,nointlimits]{article}
\usepackage{amssymb,amsmath,amsthm,verbatim}
\usepackage[dvips,final]{graphicx}
\usepackage{psfrag}
\usepackage{natbib}
\newcommand{\R}{\mathbb{R}}
\newcommand{\Prob}{\mathbb{P}}
\renewcommand{\Pr}{\Prob}
\newcommand{\supp}{\mbox{supp}}
\newcommand{\E}{\mathbb{E}}

\newcommand{\F}{\mathcal{F}}
\newcommand{\Fc}{\mathcal{F}}
\newcommand{\Lc}{\mathcal{L}}
\newcommand{\half}{\frac{1}{2}}

\newcommand{\eps}{\varepsilon}
\newtheorem{theorem}{Theorem}
\newtheorem{prop}[theorem]{Proposition}
\newtheorem{lemma}[theorem]{Lemma}
\newtheorem{cor}[theorem]{Corollary}

\theoremstyle{definition}

\newtheorem{remark}[theorem]{Remark}
\newtheorem{definition}[theorem]{Definition}

\newcommand{\indic}[1]{\boldsymbol{1}_{\{\ensuremath{#1}\}}}

\newcommand{\intsub}[1]{\int_{\{\ensuremath{#1}\}}}

\newcommand{\as}{a.s.}

\newcommand{\eg}{e.g.}

\begin{document}
\title{Skorokhod embeddings, minimality and non-centred target distributions}
\author{A.~M.~G.~Cox\thanks{e-mail:
        \texttt{amgc500@york.ac.uk}; web:
        \texttt{www-users.york.ac.uk/$\sim$amgc500/}}\\
        Department of Mathematics,\\
        University of York,\\
        York Y010 5DD, U.~K.
 \and D.~G.~Hobson\thanks{e-mail:
        \texttt{dgh@maths.bath.ac.uk}; web:
        \texttt{www.bath.ac.uk/$\sim$masdgh/}},\\
        Department of Mathematical Sciences,\\
        University of Bath,\\
        Bath BA2 7AY, U.~K.}
\maketitle

\begin{abstract}
In this paper
we consider the Skorokhod embedding problem for target
distributions with non-zero mean. In the zero-mean case, uniform
integrability provides a natural restriction on the class of
embeddings, but this is no longer suitable when the target
distribution is not centred. Instead we restrict our class of
stopping times to those which are minimal, and we find conditions
on the stopping times which are equivalent to minimality.

We then apply these results, firstly to the problem of embedding
non-centred target distributions in Brownian motion, and secondly to
embedding general target laws in a diffusion.

We construct an embedding (which reduces to the Azema-Yor embedding in the
zero-target mean case)
which maximises the law of $\sup_{s \le T} B_s$ among the class
of minimal embeddings of a general target distribution $\mu$ in
Brownian motion.
We then construct a minimal
embedding of $\mu$ in a diffusion $X$ which maximises the law of
$\sup_{s \le T} h(X_s)$ for a general function $h$.
\end{abstract}

\section{Introduction}
The Skorokhod embedding problem was first proposed, and then
solved, by \citet{Skorokhod:65}, and may be described thus:
\begin{quote}
Given a Brownian motion $(B_t)_{t \ge 0}$ and a centred target law
$\mu$ can we find a stopping time $T$ such that $B_T$ has
distribution $\mu$?
\end{quote}
Skorokhod gave an explicit construction of the stopping time $T$
in terms of independent random variables, and in doing so showed
that any zero-mean probability measure may be embedded in Brownian
motion. Since the problem was posed many more solutions have been
given, see for example \citet{Dubins:68}, \citet{Root:69} and
\citet{ChaconWalsh:76}, and the comprehensive survey article of 
\citet{Obloj:04}. With different solutions comes the
question of optimal properties of the embeddings, and various
optimal embeddings have been found --- for example the embedding
minimising the variance of $T$ \citep{Rost:76}, the embedding
minimising in convex order the law of the local time at zero
\citep{Vallois:92}, or the embedding stochastically minimising the
law of the maximum \citep{Perkins:86}.

An obvious extension of the problem is to consider more general
classes of processes. Here the question of the existence of an
embedding becomes more interesting. For a Markov process and an
arbitrary target measure necessary and sufficient conditions are
given by \citet{Rost:71} and a construction in this general case
is given by \citet{Bertoin:92}. In the case of diffusions on $\R$
simpler necessary and sufficient conditions are given in
\citet{Granditsfalkner:00}, \citet{PedersenPeskir:01} and
\citet{CoxHobson:02}, along with some constructions.

The work we present here was motivated by the following question:
\begin{quote}
Given a diffusion $(X_t)_{t \ge 0}$ and a target distribution
$\mu_X$ for which an embedding exists, which embedding maximises
the law of $\sup_{s \le T} X_s$ (respectively $\sup_{s \le T}
|X_s|$)?
\end{quote}
For Brownian motion, the question has been solved by
\citet{AzemaYor:79} (respectively \citet{Jacka:88}) under the
condition that $B_{t \wedge T}$ is a UI-martingale.

There are several considerations that need to be made when moving
from the Brownian case to the diffusion case. Firstly, the
mean-zero assumption that is made by \citet{AzemaYor:79} and
\citet{Jacka:88} is no longer natural since we are no longer
necessarily dealing with a martingale. The second aspect that
needs to be considered is with what restriction should we replace
the UI condition? That such a condition is desirable may be seen
by considering a recurrent diffusion. Here the maximisation
problem can easily seen to be degenerate by considering first
running the diffusion until it hits a level $x$, allowing it to
return to its starting point and then using the reader's favourite
embedding. Clearly this dominates the unmodified version of the
reader's favourite embedding.

In \citet{PedersenPeskir:01} an integrability condition on the
maximum (specifically that $\E(\sup_{t \le T} s(X_t)) < \infty$
where $s$ is the scale function of $X$) was suggested to replace
the UI condition in the Brownian case. In this work we propose
using the following class of stopping times introduced by
\citet{Monroe:72} to provide us with a natural restriction on the
set of admissible embeddings.
\begin{definition}
A stopping time $T$ for the process $X$ is {\em minimal} if
whenever $S \le  T$ is a stopping time such that $X_S$ and $X_T$
have the same distribution then $S = T$ \as.
\end{definition}
The class of minimal stopping times provides us with a natural
class of `good' stopping times. In the Brownian case it includes 
as a subclass those embeddings for which $\E[ \sup_{s \leq T} 
B_s ]$ is finite. Furthermore, there is a
link to the uniformly integrable case as a consequence of
the following result:
\begin{theorem}{\citep[Theorem 3]{Monroe:72}}\label{thm:monroe}
Let $T$ be an embedding in Brownian motion of a centred distribution 
$\mu$. Then $T$ is minimal if and only if the process $B_{t \wedge T}$ is 
uniformly integrable, or equivalently if $\E[B_T | {\mathcal F}_S] = B_S$
for all stopping times $S \leq T$.
\end{theorem}
Our first result extends this theorem to non-centred target distributions.
\begin{theorem}\label{thm:main1}
Let $T$ be an embedding in Brownian motion of a non-centred distribution
$\mu$ with mean $m<0$. Then $T$ is minimal if and only if the process 
$B^-_{t \wedge T}$ is
uniformly integrable, or equivalently if $\E[B_T | {\mathcal F}_S] \leq 
B_S$ for all stopping times $S \leq T$.
\end{theorem}

It is clear that the notion of a minimality fits well with the
problem of embedding in any process, such as a diffusion, and not just 
Brownian motion. Our approach to embeddings in
diffusions can be traced back to \citet{AzemaYor:79b} and will be to map 
the diffusion into natural scale (so
that, up to a time change, it resembles Brownian motion) and use
techniques developed for embedding Brownian motion. Using this
method on a transient diffusion one finds that the state space and
target distribution for the Brownian motion is restricted to a
half-line (or sometimes a finite interval). We will show
minimality to be equivalent to stopping the Brownian motion before
it leaves this interval, so that a minimal stopping time is
necessarily before the first explosion time of $X$.

When we map from the problem of embedding $\mu_X$ in $X$ to the
Brownian motion the target law $\mu$ we obtain for $B$ is the
image of $\mu_X$ under the scale function. The key point is that
there is no reason why this target law should have mean zero.
Thus, unlike most of the other studies of Skorokhod embeddings in
Brownian motion we are interested in non-centred target
distributions, and non-UI stopping times. 
Instead, as described in Theorem~\ref{thm:main1}, the class of `good' 
stopping times satisfy slightly different integrability conditions. 

Having characterised the class of minimal embeddings, we then turn to the 
problem of constructing optimal embeddings in Brownian 
motion for non-centred target distributions. In particular, amongst the 
class of minimal embeddings we find the stopping time which maximises the 
law of the maximum of the stopped process. In the centred case this 
embedding reduces to the classical Azema-Yor embedding.
In fact
most of the paper will concentrate on embedding non-centred target
distributions in $B$, and we will only return to the diffusion
case in a short final section.

The paper will proceed as follows. 
In Section~\ref{sec:minimal} we prove some
results concerning minimality of stopping times for non-centred
target distributions, giving equivalent conditions to minimality
in terms of the process. 
Next, in Section~\ref{sec:maxminembed} we construct an extension of the
Azema-Yor embedding for non-centred target distributions and show
both that it is minimal, and that it retains the optimality
properties of the original Azema-Yor embedding. In
Section~\ref{sec:maxmod} we use these stopping times to construct
an embedding maximising the distribution of $\sup_{s \le T}
h(B_s)$ for a general function $h$. Finally in
Section~\ref{sec:diffusion} we apply these results to the problem
of embedding optimally in diffusions.

\section{Minimal Embeddings for Non-centred Distributions}
\label{sec:minimal}

In this section we examine the properties of minimal stopping times. In 
particular, for embeddings in Brownian motion we aim to find equivalent 
conditions to minimality 
when the target distribution 
is not centred.

We begin by noting the following result from \citet{Monroe:72} which 
justifies the existence of minimal stopping times: 
\begin{prop}[\citet{Monroe:72}, Proposition 2] \label{prop:existenceB} For 
any stopping time $T$ there exists a minimal stopping time $S 
\le T$ 
such that $B_S \sim B_T$. \end{prop}

For $\beta \in \R$ define $H_\beta = \inf\{t > 0 : B_t = \beta\}$, the 
first hitting
time of the level $\beta$.
The main result of this section is the following:

\begin{theorem}
\label{thm:mainB} Let $T$ be a stopping time of Brownian motion
which embeds an integrable distribution $\mu$ where $m = \int_\R x \, 
\mu(dx)
< 0$. Then the following conditions are equivalent:
\begin{enumerate}
\item $T$ is minimal for $\mu$;
\item for all stopping times $R \le S \le T$,
\begin{equation}\label{eqn:Scondineq0}
\E(B_S | \Fc_R) \le B_R \text{\quad \as{};}
\end{equation}
\item for all stopping times $S \le T$,
\begin{equation}\label{eqn:Scondineq1B}
  \E (B_T | \F_S) \le B_S \text{\quad \as{};}
\end{equation}
\item for all $\gamma >0$
\begin{equation*}
\E(B_T;T>H_{-\gamma}) \le -\gamma \Pr(T>H_{-\gamma});
\end{equation*}
\item as $\gamma \to \infty$
\begin{equation*}
\gamma \Pr(T > H_{-\gamma}) \to 0;
\end{equation*}
\item
the family $\{B_S^-\}$ taken over stopping times $S \le T$ is uniformly 
integrable;
\item
for all $x > 0$
\begin{equation*}
\E(B_{T \wedge H_x}) = 0.
\end{equation*}
\end{enumerate}
In the case where $\supp(\mu) \subseteq [\alpha,\infty)$ for some
$\alpha < 0$ then the above conditions are also equivalent to the
condition:

\noindent
{\it(viii)}
\begin{equation}\label{eqn:TleqHaA}
  \Pr(T \le H_\alpha) = 1.
\end{equation}
\end{theorem}

\begin{remark} \label{rem:minus}
(i) Of course the Theorem may be restated in the case where $m > 0$ by
considering the process $-B_t$. We will use this observation
extensively in Section~\ref{sec:maxminembed}.


\noindent(ii) Equation~\eqref{eqn:Scondineq1B} is suggestive of the fact 
that when $T$ 
is minimal, 
$(B_{(u/1-u) \wedge T})_{0 \leq u \leq 1}$ is a supermartingale. To check 
this we need to 
show also that $\E B_{t \wedge T}^- < \infty$ for all $t$.
We show this more generally, for a stopping time $S \le T$. Using 
\eqref{eqn:Scondineq1B},
\begin{equation*}
\E(B_T; B_T \le 0) \le \E(B_T; B_S \le 0) \le \E(B_S; B_S \le 0),
\end{equation*}
so that $\E B_S^- < \E B_T^- < \infty$ and the process is indeed a 
supermartingale.
\end{remark}

It follows from comparing parts $(ii)$ and $(iii)$ of Theorem~\ref{thm:mainB} 
that if $S \le T$ is a stopping time and $T$ is minimal, then $S$ is 
minimal provided $\E B_S \leq 0$ and $\E |B_S| < \infty$. The first 
condition is a trivial consequence of \eqref{eqn:Scondineq0} on taking 
$R=0$, the second condition then follows from Remark~\ref{rem:minus}(ii) 
on noting that, since $\E B_S \leq 0$, we have $\E B_S^+ \leq \E B_S^- < 
\infty$.

Consequently we have the following corollary of Theorem~\ref{thm:mainB}:
\begin{cor}
If $T$ is minimal and $S \le T$ for a stopping time $S$ then $S$ is 
minimal for $\Lc(B_S)$.
\end{cor}

For ease of exposition we divide the proof of Theorem~\ref{thm:mainB} into 
a series of smaller results. The first is a key result which shows that 
the strongest of the stopping time conditions is sufficient for 
minimality. Throughout this section it is to be understood that $\mu$ is a 
distribution with negative mean and $T$ a stopping time embedding $\mu$.

\begin{lemma} \label{lem:dgh1} 
Suppose that for all stopping times $R,S$ with $R \le S \le T$ 
we have
\begin{equation} \label{eqn:Scondineq3bB}
 \E (B_S | \F_{R}) \le B_{R} \text{\quad \as{}}.
\end{equation}
Then $T$ is minimal.
\end{lemma}
                                                                                
\begin{proof}
Let $R \le T$ be a stopping time 
such that $R$ embeds $\mu$ (so that $\E|B_R|
= \E|B_T| < \infty$). For $a \in \R$,
\begin{eqnarray}
  \sup_{A \in \Fc_T} \E(a-B_T ;A)
    & = & \E(a-B_T ; B_T \le a) \nonumber\\
    & = & \E(a-B_R ; B_R \le a) \nonumber\\
    & \le &\E(a - B_T ; B_R \le a) \label{eqn:TlessSB}\\
    & \le &\sup_{A \in \Fc_T} \E(a - B_T;A) \nonumber
\end{eqnarray}
where we use \eqref{eqn:Scondineq3bB} to deduce \eqref{eqn:TlessSB}.
However since we have equality in the first and last expressions, we
must also have equality throughout and so
\begin{equation*}
 \{B_T < a \} \subseteq \{B_R \le a \} \subseteq \{B_T \le a\}.
\end{equation*}
Since this holds for all $a \in \R$ we must have $B_T = B_R$
\as{}.

It follows that for $S$ with $R \leq S \leq T$, $B_R = \E[B_R | {\mathcal 
F}_S] = \E[B_T| {\mathcal F}_S] \leq B_S$,
and since $\E[B_S | {\mathcal F}_R] \leq B_R$ we must have
$B_S=B_R=B_T$.
Thus $B$ is constant on the interval $[R,T]$ and $R=T$. Hence $T$ is
minimal.
                                                                                
\end{proof}


\begin{lemma} \label{lem:SleTexpect}
If $T$ is a stopping time such that $B_T \sim \mu$ 
and
$ 
\gamma \Pr(T > H_{-\gamma}) \to 0
$ 
as $\gamma \to \infty$ then $\E |B_S| < \infty$ and $\E B_S \le 0$ for all stopping times $S \le T$.
\end{lemma}

\begin{proof} 
We show that, for $S \le T$, $\E B_S^- < \infty$ and $\E B_S^+ \le \E B_S^-$ from which the result follows. 
Suppose $\gamma >0$. Since $B_{t \wedge H_{-\gamma}}$ is a supermartingale,
\begin{equation*}
\E( B_{T \wedge H_{-\gamma}} ; B_S < 0, S < H_{-\gamma}) \le \E(B_{S \wedge H_{-\gamma}} ; B_S < 0, S <H_{-\gamma}).
\end{equation*}
We may rewrite the term on the left of the equation as
\begin{equation*}
\E(B_T; B_S <0, T < H_{-\gamma}) - \gamma \Pr(B_S <0, S < H_{-\gamma}<T),
\end{equation*}
and by hypothesis 
$ 
\gamma \Pr(B_S <0, S < H_{-\gamma}<T) \le \gamma \Pr(H_{-\gamma} < T) \to 0
$ 
as $\gamma \to \infty$. Further, by dominated convergence,
$ 
\E(B_T; B_S < 0, T \le H_{-\gamma}) \to \E(B_T ; B_S <0)
$ 
and it follows that
\begin{equation*} 
\E (B_S ; B_S <0) = \lim_{\gamma \to \infty} \E(B_S; B_S <0, S < H_{-\gamma})\\
\ge  \E(B_T; B_S < 0).
\end{equation*}
Hence $\E B_S^- \le -\E(B_T; B_S <0) \le \E B_T^- < \infty$.

Again using the fact that $B_{t \wedge H_{-\gamma}}$ is a supermartingale,
$ 
0 \ge \E(B_{S} \wedge H_{-\gamma}) = \E(B_S ; S < H_{-\gamma}) - \gamma \Pr(H_{-\gamma} \le S)
$ 
so that
\begin{equation*}
\E(B_S^+; S < H_{-\gamma}) \le \E(B_S^- ; S < H_{-\gamma}) + \gamma \Pr(H_{-\gamma} \le S).
\end{equation*}
By monotone convergence the term on the left increases to $\E B_{S}^+$, 
while by monotone convergence and the hypothesis of the 
Lemma 
the right hand side converges to $\E B_S^-$. Consequently
$ 
\E B_S^+ \le \E B_S^- < \infty.
$ 
\end{proof}

The previous two lemmas are crucial in determining sufficient conditions 
for minimality. Now we consider necessary conditions.

\begin{lemma} \label{lem:TgeHgamma}
If $T$ is minimal then, for all $\gamma \le 0$,
$ 
  \E(B_T - B_{T \wedge H_\gamma}) \le 0.
$ 
\end{lemma}

\begin{proof}
For $\gamma \leq 0$ let $f(\gamma)= \E(B_T - B_{T \wedge H_\gamma})= 
\E(B_T - \gamma;
T > H_\gamma)$. Note that $f(0) = m < 0$. 

It is easy to see that for $0 \geq \gamma \geq \gamma'$
\begin{equation*}
  \Pr(T \in (H_{\gamma}, H_{\gamma'})) =
\Pr (\inf_{s \leq T} B_s \in (\gamma', \gamma))
\to 0
\end{equation*}
as $\gamma \downarrow \gamma'$ or $\gamma' \uparrow \gamma$.  
Since we may write
\begin{equation*}
f(\gamma) - f(\gamma') = \E(B_T; H_\gamma < T < H_{\gamma'}) + 
(\gamma' - \gamma) \Pr(T \ge H_{\gamma'}) - \gamma 
\Pr(H_{\gamma} < T < H_{\gamma'})
\end{equation*}
it follows from the
dominated convergence theorem 
that $f$ is continuous.
As a corollary if $f(\gamma_0) > 0$ for some $\gamma_0
<0$, then there exists $\gamma_1 \in (\gamma_0,0)$ such that
$f(\gamma_1)=0$.

Given this $\gamma_1$, and conditional on $T > H_{\gamma_1}$, let
$T'' = T - H_{\gamma_1}$, $W_t = B_{H_{\gamma_1} + t} -
{\gamma_1}$, and $\mu'' = \Lc(W_{T''})$.  Suppose that $T''$ is
not minimal, so there exists $S'' \le T''$ with law $\mu''$. If we
define
\begin{equation*}
  S = \begin{cases}
    T & \mbox{ on } T \le H_{\gamma_1} \\
    H_{\gamma_1} + S'' & \mbox{ on } T > H_{\gamma_1}
\end{cases} \end{equation*}
then $S$ embeds $\mu$ and $S \le T$ but $S \neq T$, contradicting
the minimality of $T$. Hence $T''$ is minimal. But then by
Theorem~\ref{thm:monroe}, $W_{t \wedge T''}$ is uniformly
integrable and so, for $\gamma < {\gamma_1}$
\begin{equation*}
  \E(W_{T''} - (\gamma -{\gamma_1});
  T'' > H_{\gamma-{\gamma_1}}^W) = 0
\end{equation*}
or equivalently
\begin{equation*}
  f(\gamma) = \E(B_T - \gamma ; T > H_{\gamma}) = 0.
\end{equation*}
Hence $f(\gamma) \leq 0$ for all $\gamma \in (-\infty,0]$.
\end{proof}

\begin{lemma} \label{lem:BTHx=0}
Suppose $B_T \sim \mu$ and $m < 0$.
If $\E[B_{T \wedge H_x}]=0$ for all $x>0$ then we have 
$\E[B_T | {\mathcal F}_S ] \leq B_S$ for all stopping times $S \leq T$.
\end{lemma}
\begin{proof}
Fix $x>0$ and define $T_x = T \wedge H_x$.
We begin by showing that $T_x$ is minimal for the centred probability
distribution ${\mathcal L}(B_{T_x})$.
Fix $R \leq T_x$. The stopped process
$B_{(t/1-t) \wedge T_x}$ is a submartingale so $\E[B_{T_x}| {\mathcal 
F}_R] \geq
B_R$ and $\E[B_R] \geq 0$. Thus 0 = $\E[B_{T_x}] \geq \E[B_R] \geq 0$ and
$B_R = \E[B_{T_x}| {\mathcal F}_R]$. Since
$B_{t \wedge T_x}$ is UI, by Theorem~\ref{thm:monroe} we have that $T_x$ 
is minimal.
 
Now fix $S \leq T$ and define $S_x = S \wedge H_x$.
We show $\E[|B_S|] \leq \E[|B_T|] < \infty$.
We have
\begin{eqnarray*}
\E[|B_S| I_{ \{ S \leq H_x \} }] & = &
\E[B_{S_x} I_{ \{ S \leq H_x \} } I_{ \{ B_S \geq 0 \} }]
- \E[B_{S_x} I_{ \{ S \leq H_x \}}  I_{ \{ B_S < 0 \} }] \\
& = & \E[B_{T_x} I_{ \{ S \leq H_x \} } I_{ \{ B_S \geq 0 \} }]
- \E[B_{T_x} I_{ \{ S \leq H_x \}}  I_{ \{ B_S < 0 \} }] \\
&\leq& \E[|B_{T_x}|  I_{ \{ T \leq H_x \} } ]
\end{eqnarray*}
Then, by two applications of Monotone Convergence
$ \E[|B_S|] \le  \E[|B_T|] < \infty$.
 
To complete the proof we show that for $A \in {\mathcal F}_S$,
$\E[ B_S I_A ] \geq \E[B_T I_A]$.
Let $A_x = A \cap \{ S \leq H_x \}$.
Then $\E[ B_S I_A ] =  \lim_x \E[ B_{S_x} I_{A_x} ] = \lim_x \E[
B_{T_x}
I_{A_x} ]$ by Monotone convergence and the minimality of $T_x$
respectively. Further,
\[
\E[ B_{T_x} I_{A_x} ]
 \geq   
\E[B_{T_x} I_{A} I_{ \{ S \leq T \leq H_x \} } ] 
 \rightarrow \E[ B_T I_A ].
\]
\end{proof}

We now turn to the proof of the main result:
\begin{proof}[Proof of Theorem~\ref{thm:mainB}]
We begin by showing the equivalence of conditions {\it(ii)} -- {\it(v)}. 
It is clear that {\it(ii)} $\implies$ {\it(iii)} $\implies$ {\it(iv)} 
$\implies$ {\it(v)}.
Hence it only remains to show that {\it(v)} $\implies$ {\it(ii)}. Suppose 
{\it(v)} holds and choose stopping times $R \le S \le T$ and $A \in \Fc_R$. Set $A_\gamma = A \cap \{R < H_{-\gamma}\}$. Since $B_{t \wedge H_{-\gamma}}$ is a supermartingale
\begin{equation}\label{eqn:SleR}
\E(B_{S \wedge H_{-\gamma}};A_\gamma) \le \E(B_{R \wedge H_{-\gamma}}; A_\gamma).
\end{equation}
By Lemma~\ref{lem:SleTexpect} $\E|B_R| < \infty$ and by dominated 
convergence 
the right hand side converges to $\E(B_R;A)$ as $\gamma \to \infty$. For the term on the left we consider
\begin{equation*}
\E(B_{S \wedge H_{-\gamma}}; A_\gamma) = \E(B_S; A\cap( S < H_{-\gamma})) 
- 
\gamma \Pr(A \cap( R<H_{-\gamma} < S)).
\end{equation*}
Again by Lemma~\ref{lem:SleTexpect} and dominated convergence the 
first term on the right converges to $\E(B_S; A)$ while the other term converges to $0$ by {\it(v)}. Hence on letting $\gamma \to \infty$ in \eqref{eqn:SleR} we have
\begin{equation*}
\E(B_S;A) \le \E(B_R;A)
\end{equation*}
and we have shown {\it(ii)}.

We have already shown that minimality is equivalent to these conditions: 
{\it(ii)} $\implies$ {\it(i)} is Lemma~\ref{lem:dgh1}, 
while {\it(i)} $\implies$ {\it(iv)} is Lemma~\ref{lem:TgeHgamma}.

Now consider {\it(vi)}.
If \eqref{eqn:Scondineq1B} holds then $B_S^- \leq \E(B_T|{\mathcal F}_S)^-
\leq \E(B_T^- | {\mathcal F}_S )$ and uniform integrability follows.
Conversely,
{\it(vi)} implies
$\sup_{S \le T} \E(B_S^-; B_S^- \ge \gamma) \to 0$ as $\gamma \to \infty$. 
Taking $S = H_{-\gamma} \wedge T$ yields {\it(v)}.

The implication {\it(vii)} implies {\it(iii)} is the content of 
Lemma~\ref{lem:BTHx=0}. We now show {\it(vi)} implies {\it(vii)}. Fix 
$x >0$. Since 
$B_S^-$ is uniformly integrable for all $S$, we have that $B_{t \wedge H_x 
\wedge T}^-$ is uniformly integrable. Since $B_{t \wedge H_x
\wedge T}$ is also bounded above, it follows that $B_{t \wedge H_x
\wedge T}$ is UI. Hence $\E[B_{H_x \wedge T}] = \lim_t \E[ B_{t \wedge H_x 
\wedge T}]=0$. 

We have shown equivalence between {\it (i)} -- {\it (vii)}. 
We are left with showing that if $\mu$ has
support bounded below then {\it (viii)} is also equivalent. So assume
that the target distribution $\mu$ has support contained in
$[\alpha, \infty)$ for some $\alpha < 0$ and that $T$ is an embedding of 
$\mu$. 
Then $B_{t \wedge H_\alpha}$ is a continuous supermartingale,
bounded below and therefore if $S \le T \le H_\alpha$,
\begin{equation*}
  \E(B_T|\Fc_S) \le B_S.
\end{equation*}
The reverse implication follows from considering the stopping time
$H_{\alpha - \epsilon} = \inf\{t \ge 0 : B_t \le \alpha - \eps\}$,
for then if $A= \{ \omega : H_{\alpha - \epsilon} < T \}$ and $S =
H_{\alpha  - \eps} \wedge T$,
\begin{equation*}
  (\alpha - \epsilon) \Pr(A) = \E(B_{H_{\alpha - \epsilon}};A) =
  \E(B_S;A) \geq \E(B_T ;A) \geq \alpha \Pr(A)
\end{equation*}
which is only possible if $\Pr(A)=0$.
\end{proof}


We close this section with a discussion of the case where the mean of the 
target distribution is well defined, but infinite. In this case the notion 
of minimality still makes perfect sense, but many of the 
conditions outlined in Theorem~\ref{thm:mainB} can be shown to be no 
longer equivalent to minimality.

Suppose that $\mu$ only places mass on $\R^-$ and that 
that $\mu( (-\infty,x)) > (1 + |x|)^{-1}$ for $x < 0$. 
Let $\phi : \R^- \mapsto \R^-$ solve  
\[ \mu( (-\infty,x) ) = \frac{1}{1 + |\phi^{-1}(x)| }. \]
Such a $\phi$ is an increasing function with $\phi(x) < x$.

Define $J = \inf \{B_u ; u \leq H_1 \}$. By construction $\phi(J) \sim 
\mu$ and we can construct an embedding of $\mu$ by setting $T = \inf \{u > 
H_1 : B_u = \phi(J) \}$. Further, at the stopping time $T$ the 
Brownian motion is at a current minimum. Hence $T$ is minimal.

Now consider the alternative conditions in Theorem~\ref{thm:mainB}. In 
reverse order, {\it(vii)} and {\it(vi)} both clearly fail, whereas 
{\it(v)} may or may not hold depending on the target law $\mu$, and 
{\it(iv)} holds trivially. (However {\it(iv)} holds for all embeddings of $\mu$
so that it is not sufficient for minimality.) The condition in 
${\it(iii)}$ also holds, 
since on stopping $B$ is at a minimum value. However the choice $R=0$ and 
$S=H_1$ shows that ${\it(ii)}$ fails and is not necessary for minimality.
Hence, in the case where $m = \int_{\R} x \mu(dx)$ is well defined but 
equal to $- \infty$, the only condition from the list in 
Theorem~\ref{thm:mainB} which could possibly 
be necessary and sufficient for minimality is ${\it(iii)}$.

\section{An extension of the Az\'{e}ma-Yor embedding to the non-centred 
case} \label{sec:maxminembed}

Let $\mu$ be a target distribution on $\R$, and let $B_t$ be a Brownian 
motion with $B_0=0$. Define $M_T = \sup_{s \leq T}B_s$ and $J_T= \inf_{s 
\leq T}B_s$. In this section we are interested in finding an embedding to 
solve the following problem:
\begin{quote}
  Given a Brownian motion $(B_t)_{t \ge 0}$ and an integrable (but possibly
  not centred) target distribution $\mu$ with mean $m$, find a minimal
  stopping time $T$ such that $T$ embeds $\mu$ and
  \begin{equation*}
    \Pr(M_T \ge x)
  \end{equation*}
  is maximised over all minimal stopping times $T$ embedding
  $\mu$, uniformly in $x$.
\end{quote}
We call an embedding with this property the max-max embedding, and
denote it by $T_{max}$.

Without some condition on the class of admissible stopping times
the problem is clearly degenerate --- any stopping time may be
improved upon by waiting for the first return of the process to 0
after hitting level $x$ and then using the original embedding. For
this improved embedding $\Pr (M_T \geq x) = 1$. Further, since no
almost surely finite stopping time can satisfy
$ 
  \Pr(M_T \ge x) = 1
$ 
for all $x >0$, there can be no solution to the problem above in
the class of all embeddings. As a consequence some restriction on
the class of admissible stopping times is necessary for us to have
a well defined problem.

Various conditions have been proposed in the literature to
restrict the class of stopping times. In the case where $m=0$, the
condition on $T$ that $B_{t\wedge T}$ is a UI martingale was
suggested by \citet{DubinsGilat:78}, and in this case the maximal
embedding is the Azema-Yor embedding. When $m=0$ \citet{Monroe:72}
tells us that minimality and uniform integrability are equivalent
conditions, so the Azema-Yor stopping time is the max-max
embedding. For the case where $m>0$, \citet{PedersenPeskir:01}
showed that $\E M_T < \infty$ is another suitable condition, with
the optimal embedding being based on that of Azema and Yor. We
argue that the class of minimal embeddings is the appropriate
class for the problem under consideration since minimality is a
natural and meaningful condition, which makes sense for all $m$
(and which, for $m>0$, includes as a subclass those embeddings
with $\E(M_T) < \infty$.)

\psfrag{a}{$\alpha$}
\psfrag{b}{$\beta$}
\psfrag{c(x)}{$c(x)$}
\psfrag{x}{$x$}
\psfrag{z}{$z$}
\psfrag{b(z)}{$b(z)$}
\psfrag{2m}{$2m$}
\psfrag{z-}{$-z_-(\theta_0)$}
\psfrag{z+}{$z_+(\theta_0)$}
\psfrag{u}{$u(\theta_0)$}
\begin{figure}[t]
\begin{center}
\includegraphics[width=\textwidth,height=3in]{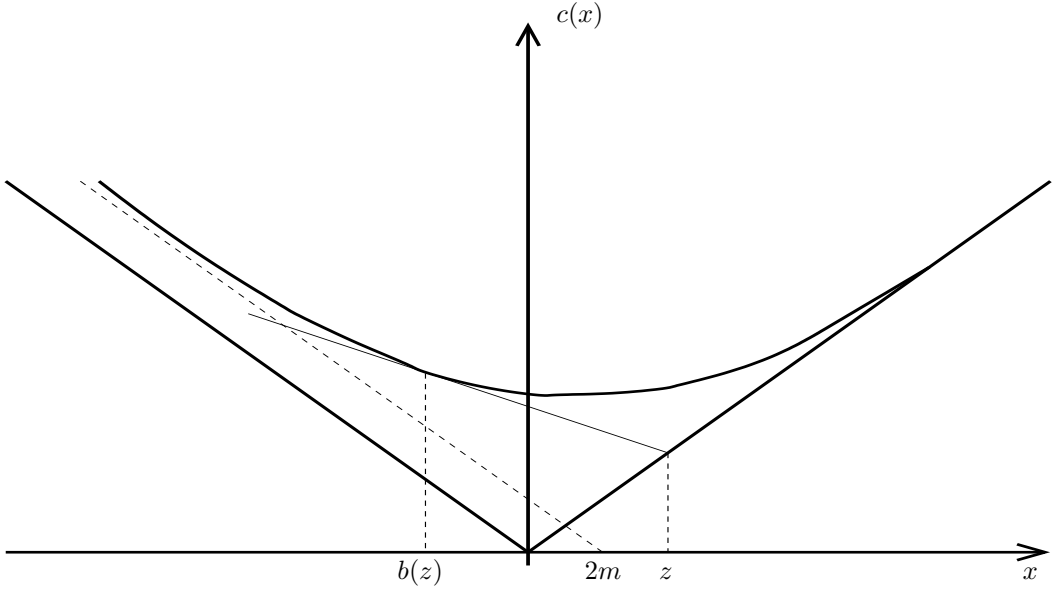}
\caption{\label{fig:cgraph2} $c(x)$ for a $\mu$ with support
bounded above, and positive non-zero mean $m$. Also shown is an
intuitive idea of $b(z)$.}
\end{center}
\end{figure}


We now describe the construction of the candidate max-max stopping
time. This construction is an extension of the classical 
\citet{AzemaYor:79} stopping time. However we derive the stopping rule 
in a slightly non-standard way via a \citet{ChaconWalsh:76} style 
argument. This interpretation of the Az\'{e}ma-Yor construction (in 
the 
centred case) is due to \citet{Meilijson:83}.

In the following we treat the cases $m>0$, $m=0$ and $m<0$ in one go, 
since the basic idea is identical.
Define the convex function:
\begin{equation} \label{eqn:cdefn2}
  c(x) := \E^\mu |X-x| +|m|.
\end{equation}
We note that $c$ is related to the potential of $\mu$, that it is 
Lebesgue-almost everywhere differentiable with left-derivative 
$c'_-(x) = 1 - 2 \mu([x, \infty))$ and
$c(x)-|x| \to |m| \mp m$ as $x\to \pm \infty$. 
For $\theta \in [-1,1]$ define
\begin{equation*}
  u(\theta) := \inf \{ y \in \R : c(y) + \theta(x-y) \le c(x),
    \forall x \in \R \},
\end{equation*}
\begin{equation*}
  z_+(\theta) := \frac{c(u(\theta)) - \theta u(\theta)}{1-\theta},
\end{equation*}
and (note that $z^{-1}_+$ is well defined) for $x \ge 0$
\begin{equation*}
  b(x) := u(z_{+}^{-1}(x)).
\end{equation*}
The intuition behind these definitions 
is as follows. Let $\theta$ denote a gradient and consider the unique 
tangent to $c$ with gradient $\theta$. Then $z_+(\theta)$ is the 
$x$-coordinate of the point where this tangent crosses the line $y=x$.
Similarly, $b(x)$ is the $x$-coordinate 
of the point for which the tangent to $c$ at that point passes through 
$(x,x)$.
The mathematical definitions above, and the fact that $u$ is 
left-continuous, ensure 
that $b$ is well defined. Note that when $m=0$, a few lines of 
calculus are sufficient to check that $b^{-1}$ is precisely the 
barycentre function which arises in the classical Azema-Yor construction.

\begin{theorem} \label{thm:AYmNonZero}
Let $T$ be a stopping time of $(B_t)_{t \ge 0}$ which embeds $\mu$
and is minimal. Then for $x \ge 0$
\begin{equation} \label{eqn:Jineq}
  \Pr(M_T \ge x) \le \half \inf_{\lambda < x}
    \left( \frac{c(\lambda) - \lambda}{x - \lambda}\right).
\end{equation}
Define the stopping
time $T_{max}$ via
\begin{equation} \label{eqn:TMMdefn}
  T_{max} := \inf\{t > 0 : B_t \le b(M_t) \}.
\end{equation}
Then $T_{max}$ embeds $\mu$, is 
minimal, and attains 
equality in
\eqref{eqn:Jineq} for all $x \ge 0$.
\end{theorem}

\begin{remark} (i) By the comments before the theorem relating $b$ to the 
barycentre function, when $m=0$ the above theorem is a restatement of the 
classical Azema-Yor result, which by Theorem~\ref{thm:monroe} can be 
stated in terms of minimal, rather than uniformly integrable embeddings. 
For $m>0$, $b(x) = -\infty$ for $x < m$, and consequently $T_{max} \ge 
H_m$. Moreover, for $x>m$, $b^{-1}$ is a shifted version of the barycentre 
function. Consequently, when $m>0$ the embedding $T_{max}$ may be thought 
of as `wait until the process hits $m$ then use the Azema-Yor embedding,' 
and the conclusion of Theorem~\ref{thm:AYmNonZero} is similar to 
Proposition 3.1 of 
\citet{PedersenPeskir:01}, except that $T_{max}$ is shown to be optimal 
amongst the larger class of minimal stopping times rather than the class 
of embeddings for which $M_T \in L^1$. However, the truly original part of 
the theorem is in the case $m<0$. In that case the embedding `wait until 
the process hits $m$ and then use the Azema-Yor embedding' does not 
achieve equality in \eqref{eqn:Jineq}.

\noindent(ii) Note that
\begin{equation} \label{eqn:cLambda}
  \frac{c(\lambda) - \lambda}{x - \lambda}
    = 1 - \frac{x - c(\lambda)}{x - \lambda}.
\end{equation}
We can relate the right-hand-side of \eqref{eqn:cLambda} to the
slope of a line joining $(x,x)$ with $(\lambda,c(\lambda))$. In
taking the infimum over $\lambda$ we get a tangent to $c$ and a
value for the slope in $[-1,1]$. Thus the bound on the
right-hand-side of \eqref{eqn:Jineq} lies in $[0,1]$.

\noindent(iii) $T_{max}$ has the property that it maximises the law of 
$M_T$ over
minimal stopping times which embed $\mu$. If we want to minimise
the law of the minimum, or equivalently we wish to maximise the
law of $-J_{T}$, then we can deduce the form of the optimal
stopping time by reflecting the problem about 0, or in other words
by considering $-B$. Let $T_{min}$ be the embedding which arises
in this way, so that amongst the class of minimal stopping times
which embed $\mu$, the stopping time $T_{min}$ maximises
\begin{equation*}
  \Pr ( - J_T \ge x )
\end{equation*}
simultaneously for all $x \ge 0$.
\end{remark}

We now turn to the proof of Theorem~\ref{thm:AYmNonZero}.

\begin{proof}
The following inequality for $x>0$, $\lambda < x$ may be verified
on a case by case basis:
\begin{equation} \label{eqn:Jident}
  \indic{M_T \ge x} \le \frac{1}{x-\lambda}
    \left[ B_{T \wedge H_{x}}
    + \frac{|B_T - \lambda| - (B_T + \lambda)}{2}\right].
\end{equation}
In particular, on $\{M_T < x\}$, \eqref{eqn:Jident} reduces to
\begin{equation} \label{eqn:Jident2}
  0 \le \
  \begin{cases}
    0 & \lambda \ge B_T \\
    \frac{B_T - \lambda}{x- \lambda} & \lambda < B_T,
  \end{cases}
\end{equation}
and on $\{M_T \ge x\}$ we get
\begin{equation} \label{eqn:Jident3}
  1 \le
  \begin{cases}
    \frac{x- B_T}{x - \lambda} & \lambda > B_T \\
    1 & \lambda \le B_T.
  \end{cases}
\end{equation}
Then taking expectations,
\begin{equation} \label{eqn:Twedgeineq0}
  \Pr(M_T \ge x) \le \frac{1}{x-\lambda} \left[ \E (B_{T \wedge H_{x}})
    + \frac{c(\lambda) - |m| - (m + \lambda)}{2}\right].
\end{equation}
If $m\le 0$ then by Theorem~\ref{thm:mainB}{\it (vii)} and the minimality of
$T$ we know $\E(B_{T \wedge H_x}) = 0$ and so
\begin{equation*}
  \Pr(M_T \ge x) \le \half \frac{c(\lambda) - \lambda}{x - \lambda}.
\end{equation*}
Conversely if $m > 0$, by Theorem~\ref{thm:mainB}$(iii)$ applied to $-B$,
\begin{equation} \label{eqn:Twedgeineq}
  m = \E(B_T) \ge \E(B_{T \wedge H_x})
\end{equation}
and so
\begin{equation*}
  \Pr(M_T \ge x) \le \frac{1}{x-\lambda} \left[ m + \frac{c(\lambda)
    - 2m - \lambda}{2}\right] = \half
    \frac{c(\lambda) - \lambda}{x - \lambda}.
\end{equation*}
Since $\lambda$ was arbitrary in either case, \eqref{eqn:Jineq}
must hold. 

It remains to show that $T_{max}$ attains equality in
\eqref{eqn:Jineq}, embeds $\mu$ and is minimal.
We begin by showing that it does attain equality in
\eqref{eqn:Jineq}. Since
\begin{equation*}
  \frac{c(\lambda) - \lambda}{x-\lambda} =
    1 + \frac{c(\lambda)-x}{x-\lambda}
\end{equation*}
the infimum in \eqref{eqn:Jineq} is attained by a value
$\lambda^*$ with the property that a tangent of $c$ at $\lambda^*$
intersects the line $y=x$ at $(x,x)$.  By the definition of $b$ we
can choose $\lambda^* = b(x)$. In particular, since $\{M_{T_{max}}
< x\} \subseteq \{B_{T_{max}} \le b(x)\}$ and $\{M_{T_{max}} \ge
x\} \subseteq \{B_{T_{max}} \ge b(x)\}$, the stopping time
$T_{max}$ attains equality almost surely in \eqref{eqn:Jident2}
and \eqref{eqn:Jident3}. Assuming that $T_{max}$ is minimal, we
are then done for $m \le 0$. If $m>0$ we do not always have
equality in \eqref{eqn:Twedgeineq}.  If $x < m$ then
$\E(B_{T_{max} \wedge H_x}) = x$, but then $\lambda^* = -\infty$
and so the term $\E (B_{T_{max} \wedge H_x}) - m)/(x -
\lambda^*)$ in \eqref{eqn:Twedgeineq0} is zero. As a result equality is 
again attained in
\eqref{eqn:Jineq}. Otherwise, if $x \ge m$ then $T_{max} \geq H_m$
and the properties of the classical Azema-Yor embedding ensure
that $\E (B_{T_{max} \wedge H_x}) = m$ and there is equality both
in \eqref{eqn:Twedgeineq} and \eqref{eqn:Jineq}.

Fix a value of $y$ which is less than the supremum of the support
of $\mu$. Let $b^{-1}$ be the left-continuous inverse of $b$. Then,  
since we
have equality in \eqref{eqn:Jineq}, we deduce:
\begin{eqnarray*}
  \Pr(B_{T_{max}} \ge y)
    &=& \Pr (M_{T_{max}} \ge b^{-1}(y)) \\
    &=& \half \left[ 1+\frac{c(b(b^{-1}(y))) -b^{-1}(y)}{b^{-1}(y) -
      b(b^{-1}(y))} \right]\\
    &=&\half (1-c_{-}'(y))
\end{eqnarray*}
where $c'_-$ is the left-derivative of $c$. It is easy to see from the 
definition \eqref{eqn:cdefn2} that this last expression equals 
$\mu([y,\infty))$.
Hence $T_{max}$ embeds $\mu$.

Now we consider minimality of $T_{max}$. It is well known that in the 
centred case $T_{max}$ is uniformly integrable, and hence by 
Theorem~\ref{thm:monroe}, $T_{max}$ is minimal. Suppose $m>0$. By an 
analogue of Theorem~\ref{thm:mainB}$(v)$, in order prove minimality it is 
sufficient to show that $\lim_{x \uparrow \infty} x \Pr( T_{max} > H_x ) 
\rightarrow 0$. But, by the arguments in the previous paragraph 
\[ 2 x \Pr( T_{max} > H_x ) = 2 x \Pr( M_{T_{max}} > x ) = x 
\left[ 1 - c'_- (b(x)) \right] \]
This last quantity is exactly the height above 0, when it crosses the 
$y$-axis, of the tangent to $c$ which passes through $(x,x)$. As $x 
\uparrow \infty$ this height decreases to zero. 

Now suppose $m<0$. By Theorem~\ref{thm:mainB}{\it(v)}, in order to 
prove minimality it is sufficient to 
show that as $x 
\downarrow - \infty$, $|x| \Pr({T_{max}} > H_x) \rightarrow 0$. We 
have
\[ |x| \Pr( {T_{max}} > H_x ) = |x| \Pr( H_x < H_{b^{-1}(x)}) = 
\frac{|x| 
b^{-1}(x)}{|x| + b^{-1}(x)} < b^{-1}(x). \]
It is easy to see from the representation of $b$ that $b^{-1}(x)$ tends to 
zero as $x \rightarrow - \infty$. 
\end{proof}

\section{An embedding to maximise the modulus}\label{sec:maxmod}

\citet{Jacka:88} shows how to embed a centred
probability distribution in a Brownian motion so as to maximise
$\Pr (\sup_{t \leq T} |B_t| \geq y)$. Our goal in this section is
to extend this result to allow for non-centered target
distributions with mean $m \neq 0$. In fact we solve a slightly
more general problem. Let $h$ be a measurable function; we will
construct a stopping time $T_{mod}$ which will maximise
$\Pr(\sup_{t \le T} |h(B_t)| \ge y)$ simultaneously for all $y$
where the maximum is taken over the class of all minimal stopping
times which embed $\mu$.  The reason for our generalisation will
become apparent in the application in the next section.

Without loss of generality we may assume that $h$ is a
non-negative function with $h(0)=0$ and such that for $x > 0$ both
$h(x)$ and $h(-x)$ are increasing. To see this, observe that for
arbitrary  $h$ we can define the function
\begin{equation*}
  \tilde{h}(x) =
  \begin{cases}
    \max_{0 \le y \le x} |h(y)| - |h(0)| & x \ge 0; \\
    \max_{x \le y \le 0} |h(y)| - |h(0)| & x < 0.
  \end{cases}
\end{equation*}
Then $\tilde{h}$ has the desired properties and since
\begin{equation*}
  \sup_{s \le T} |h(B_s)| = \sup_{s \le T} \tilde{h}(B_s)
    + |h(0)|
\end{equation*}
the optimal embedding for $\tilde{h}$ will be an optimal
embedding for ${h}$.

So suppose that $h$ has the properties listed above.  We want to
find an embedding of $\mu$ in $B$ which is minimal and which
maximises the law of $\sup_{t \le T} h(B_t)$.  (Since $h$ is
non-negative we can drop the modulus signs.) Suppose also for
definiteness that $\mu$ has a finite, positive mean $m = \int_\R x
\mu(dx) > 0$. In fact our construction will also be optimal when
$m=0$ (the case covered by \citet{Jacka:88}), but in order to
avoid having to give special proofs for this case we will omit it.

\begin{figure}[t]
\begin{center}
\includegraphics[width=\textwidth,height=3in]{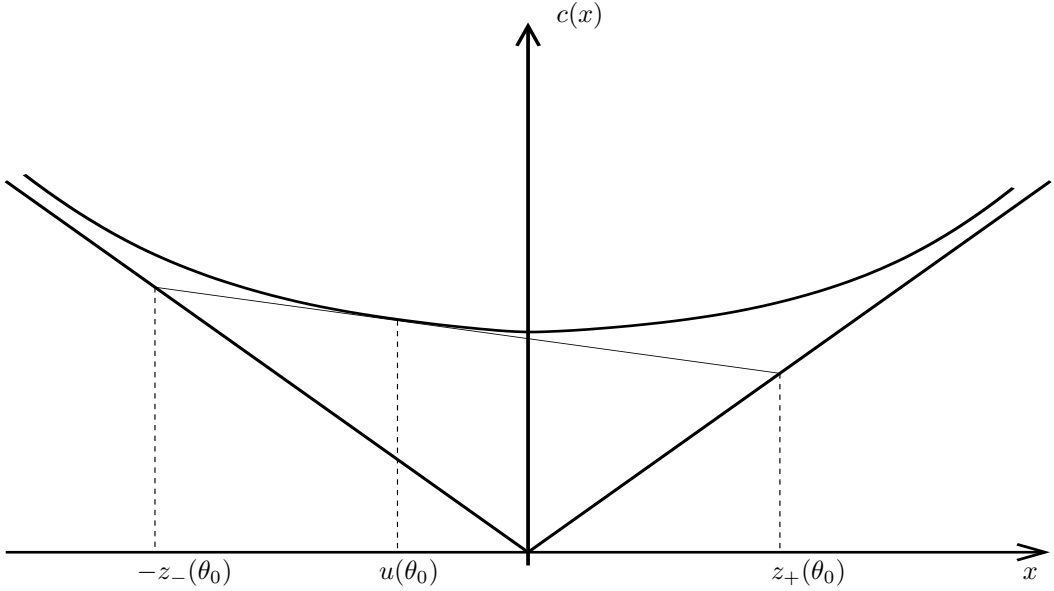}
\caption{\label{fig:cgraph4} $c(x)$ for a distribution $\mu$
showing the construction of $z_+(\theta_0)$ and $z_-(\theta_0)$.
The slope of the tangent is $\theta_0$ where $\theta_0$ has been
chosen such that (assuming $h$ is continuous) $h(z_+(\theta_0)) =
h(-z_-(\theta_0))$.}
\end{center}
\end{figure}

Recall the definitions of $c$ and $u$ from the previous section. 
Define
\begin{eqnarray}
  z_+(\theta) &:=& \frac{c(u(\theta)) - \theta u(\theta)}{1
    - \theta}, \label{eqn:zpdefn2}\\
  z_-(\theta) &:=& \frac{c(u(\theta)) - \theta u(\theta)}{1
    + \theta},\label{eqn:zmdefn2}
\end{eqnarray}
so that $-z_{(\theta)}$ is the $x$-coordinate of the point where the 
tangent 
to $c$ with slope $\theta$ intersects the line $y=-x$,
and set
\begin{equation*}
  \theta_0 := \inf \{ \theta \in [-1,1] : h(z_+(\theta))
    \ge h(-z_-(\theta))\},
\end{equation*}
as pictured in Figure~\ref{fig:cgraph4}. Our optimal stopping time
will take the following form. Run the process until it hits either
$z_+(\theta_0)$ or $-z_-(\theta_0)$, and then embed the
restriction of $\mu$ to $[u(\theta_0),\infty)$ or
$(-\infty,u(\theta_0)]$ respectively (defining the target measures
more carefully when there is an atom at $u(\theta_0)$). For the
embeddings in the second part, we will use the constructions
described in Section~\ref{sec:maxminembed}.

To be more precise about the measures we embed in the second step,
define
\begin{equation*}
  p:= \Pr(H_{z_+(\theta_0)} < H_{-z_{-}(\theta_0)}) =
    \frac{z_-(\theta_0)}{z_+(\theta_0)+z_-(\theta_0)},
\end{equation*}
and note
\begin{equation*}
  \theta_0=\frac{z_+(\theta_0)- z_-(\theta_0)}{z_+(\theta_0)
    + z_-(\theta_0)} = 1-2p.
\end{equation*}
Let $\mu_+$ (respectively $\mu_-$) be the measure obtained by
conditioning a random variable with law $\mu$  to lie in the upper
$p^{\mbox{th}}$ (respectively lower $(1-p)^{\mbox{th}}$) quantile
of its distribution.

Recall that we have taken $m>0$ and that
\begin{eqnarray}
  c(y)  = \int |w-y| \mu(dw) + |m| 
& = & 2 \intsub{w > y} (w-y) \mu(dw) + y  \label{eqn:c+def} \\
& = & 2 \intsub{w < y} (y-w) \mu(dw) + 2m - y \label{eqn:c-def}
\end{eqnarray}
Then using \eqref{eqn:c+def} in the definition \eqref{eqn:zpdefn2} we have 
that
\begin{eqnarray*}
  z_+(\theta_0) & = &
    \frac{1}{2p} \left( 2 \intsub{w > u(\theta_0)}
    (w-u(\theta_0))\, \mu(dx) + 2p u(\theta_0) \right) \\
  & = & \frac{1}{p} \intsub{w > u(\theta_0)} w \, \mu(dw) +
    u(\theta_0)\left( 1- \frac{1}{p} \mu((u(\theta_0),\infty))\right)
\end{eqnarray*}
In particular $z_+(\theta_0)$ is the
mean of $\mu_+$, since $\mu_+(\{u(\theta_0)\}) = 1-
\frac{1}{p}\mu((u(\theta_0),\infty))$. When we repeat the
calculation for $z_-(\theta_0)$ using \eqref{eqn:c-def} we find that
\begin{equation*}
  -z_-(\theta_0) = \frac{1}{1-p} \intsub{w < u(\theta_0)} w \,
    \mu(dw) + u(\theta_0) \left( 1- \frac{1}{1-p}
    \mu((-\infty,u(\theta_0)) \right) - \frac{m}{(1-p)}.
\end{equation*}
Note that $-z_-(\theta_0)$
is strictly smaller than the mean of $\mu_-$.

We now describe the candidate stopping time $T_{mod}\equiv
T^h_{mod}$. Note that this stopping time will depend implicitly on
the function $h$ via $z_{\pm}(\theta_0)$. Let
\begin{equation*}
  T_0 := \inf \{t>0:B_t \not\in (-z_{-}(\theta_0),z_+(\theta_0))\},
\end{equation*}
and define
\begin{equation*}
  T_{mod} :=
  \begin{cases}
    T_{max}^{\mu_+} \circ \theta_{T_0}
      & \quad B_{T_0} = z_+(\theta_0) \\
    T_{min}^{\mu_-} \circ \theta_{T_0}
      & \quad B_{T_0} = -z_-(\theta_0).
  \end{cases}
\end{equation*}
Here we use $\theta_{T_0}$ to denote the shift operator, and
$T_{max}^{\mu_+}$ is the stopping time constructed in
Section~\ref{sec:maxminembed} for a zero-mean target distribution,
so that $T_{max}^{\mu_+}$ is a standard Azema-Yor embedding of the
centred target law ${\mu_+}$. (Recall that $z_+(\theta_0)$ is the
mean of the corresponding part of the target distribution.)
Similarly $T_{min}^{\mu_-}$ is the stopping time applied to $-B$
started at $-z_-(\theta_0)$ which maximises the law of the maximum
of $-B$. In this case the mean of the target law $\mu_-$ is larger
than $-z_-(\theta_0)$ so that in order to define $T_{min}^{\mu_-}$
we need to use the full content of Section~\ref{sec:maxminembed}
for embeddings of non-centred distributions.

The following theorem asserts that this embedding is indeed an
embedding of $\mu$, that it is minimal, and that it has the
claimed optimality property.

\begin{theorem} \label{thm:maxmod}
Let $\mu$ be a target distribution such that $m >0$. Then within
the class of minimal embeddings of $\mu$ in $B$, the embedding
$T_{mod}$ as defined above has the property that it maximises
\begin{equation*}
  \Pr\left(\sup_{t \le T} h(B_t) \ge x\right)
\end{equation*}
simultaneously for all $x$.
\end{theorem}
\begin{proof}
By construction $T_{mod}$ embeds $\mu$. We need only show that it
is optimal and minimal.

For $x \le h(-z_-(\theta_0)) \wedge h(z_+(\theta_0))$ we know
the probability that $\{\sup_{t \le T_{mod}} h(B_t) \ge
x\}$ is one and so, for such $x$, $T_{mod}$ is clearly optimal.
Indeed if $h$ is discontinuous at $-z_-(\theta_0)$ or
$z_+(\theta_0)$ slightly more can be said. Note first that if
$z_+(\theta_0)$ coincides with the supremum of the support of
$\mu$, then by Theorem~\ref{thm:mainB}$(iii)$ and the minimality of
$T_{mod}$ (see below), the stopped Brownian motion can never go
above $z_+(\theta_0)$. Hence we may assume that $h$ is constant on the 
interval to the right of the supremum of its support.

Define
\begin{equation*}
  L = \left( \lim_{y \uparrow -z_-(\theta_0)} h(y) \right)
    \wedge \left(\lim_{y \downarrow z_+(\theta_0)} h(y) \right)
\end{equation*}
and take $x \leq L$. Then either $B_{T_0} =
z_+(\theta_0)$ or $B_{T_0} = -z_-(\theta_0)$. If $B_{T_0} =
z_+(\theta_0)$ then either $M_{T_{mod}} > z_+(\theta_0)$ almost
surely and
\begin{equation*}
  \max_{0 \le t \le T_{mod}} h(B_t) \ge \lim_{y
    \downarrow z_+(\theta_0)} h(y) \ge L
\end{equation*}
or $z_+(\theta_0)$ is the supremum of the support of $\mu$ and
\begin{equation*}
  \max_{0 \le t \le T_{mod}} h(B_t) = h(B_{T_0}) \ge L.
\end{equation*}
Similar considerations apply for $B_{T_0} = -z_-(\theta_0)$ except
that then $-J_{T_{mod}} > z_-(\theta_0)$ in all cases. We deduce
that for $x \leq L$
\begin{equation*}
  \Pr\left(\sup_{t \le T_{mod}} h(B_t) \ge x\right) = 1
\end{equation*}
and hence $T_{mod}$ is optimal for such $x$.

So suppose that $x >L$. For any stopping time $T$ embedding $\mu$,
the following holds:
\begin{equation} \label{eqn:disjprob}
  \Pr\left(\sup_{s \le T} h(B_s) \ge x\right) \le
  \Pr\left(h(M_T) \ge x \right) + \Pr\left(h(J_T) \ge x\right).
\end{equation}
We will show that the embedding $T_{mod}$ attains the maximal
values of both terms on the right hand side, and further that for
$T_{mod}$ the two events on the right hand side are disjoint.
Hence $T_{mod}$ is optimal.

By the definition of $\theta_0$, $x > \left( h(z_+(\theta_0))
\right) \vee \left( h(-z_-(\theta_0)) \right)$. It follows that
\begin{equation*}
  \Pr(h(M_{T_{mod}}) \ge x)  =
    p \Pr ( h(M_{T_{mod}}) \ge x | B_{T_0} = z_+(\theta_0))
\end{equation*}
and by the definition of $T_{mod}$ and the properties of
$T_{max}$, we deduce
\begin{eqnarray*}
  \Pr(h(M_{T_{mod}}) \ge x)
    & = & p\Pr (h(M_{T^{\mu}_{max}}) \ge x | M_{T^{\mu}_{max}}
      \ge z_+(\theta_0))\\
    & = & \Pr(h(M_{T^{\mu}_{max}}) \ge x)
\end{eqnarray*}
where here $T^{\mu}_{max}$ is the embedding of
Section~\ref{sec:maxminembed} applied to $\mu$. A similar
calculation can be done for the minimum. In particular $T_{mod}$
inherits its optimality property from the optimality of its
constituent parts $T^{\mu_+}_{max}$ and $T^{\mu_-}_{min}$

Finally we note that $T_{mod}$ is indeed minimal. 
Consider the family of stopping times $S \leq T_{mod}$. On the set where 
$B_{T_0} = - z_-(\theta_0)$ we have that $B_S^+$ 
is bounded, whereas on the set $B_{T_0} = -z_-(\theta_0)$ 
the minimality of $T^{\mu_-}_{min}$ ensures that $B_S^+$
is uniformly integrable. Hence, combining these two cases, $B_S^+$ is a 
uniformly integrable family 
and the minimality of $T_{mod}$ follows from \ref{thm:mainB}{\it(vi)} 
applied to a target distribution with positive mean. 
\end{proof}

\begin{remark}
If the restrictions of $h$ to $\R_{+}$ and $\R_{-}$ are strictly
increasing then $T_{mod}$ will be essentially the unique embedding
which attains optimality in Theorem~\ref{thm:maxmod}. If however
$h$ has intervals of constancy then other embeddings may also
maximise the law of $\sup_{t \leq T} |h(B_t)|$.
\end{remark}

\section{Embeddings in diffusions} \label{sec:diffusion}

Our original motivation in considering the embeddings of the
previous sections was their use in the investigation of the
following question:
\begin{quote}
  Given a regular (time-homogeneous) diffusion $(X_t)_{t \ge 0}$ and
  a target distribution $\mu_X$, find (if possible) a minimal
  stopping time which embeds $\mu_X$ and which maximises the law of
  $\sup_{t \le T} X_t$ (alternatively $\sup_{t \le T} |X_t|$) among
  all such stopping times.
\end{quote}

Note that in the martingale (or Brownian) case it is natural to
consider centred target laws, at least in the first instance.
However in the non-martingale case this restriction is no longer
natural, and as we shall see below is completely unrelated to
whether it is possible to embed the target law in the diffusion
$X$. It was this observation which led us to consider the problem
of embedding non-centred distributions in $B$.

The key idea (see \citet{AzemaYor:79b}) is that we can relate the problem 
of embedding in a
diffusion to the case where we are dealing with Brownian motion
via the scale function. There exists a continuous, increasing
function $s$ such that $s(X_t)$ is a local martingale, and hence a
time-change of Brownian motion.  Then the requirement $X_{T^X}
\sim \mu_X$ translates to finding an embedding of a related law in
a Brownian motion $B$, and the criterion of maximising $\sup_{t
\le T} X_t$ also has an equivalent statement in terms of $B$.

We first recall the properties of the scale function (see \eg{}
\citet[V.45]{RogersWilliams:00b}). If $(X_t)_{t \ge 0}$ is a
regular (time-homogeneous) diffusion on an interval $I \subseteq
\R$ with absorbing or inaccessible endpoints and vanishing at
zero, then there exists a continuous, strictly increasing scale
function $s:I \to \R$ such that $Y_t = s(X_t)$ is a diffusion in
natural scale on $s(I)$. We may also choose $s$ such that
$s(0)=0$. In particular $Y_t$ is (up to exit from the interior of
$s(I)$) a time change of a Brownian motion with strictly positive
speed measure. For definiteness we write $Y_t = B_{\tau_t}$.

We suppose also that our target distribution $\mu_X$ is
concentrated on the interior of $I$. Then we may define a measure
$\mu = \mu_Y$ on $s(I)^\circ$ by:
\begin{equation*}
  \mu(A) = \mu_Y(A) = \mu_X(s^{-1}(A)),
    \mbox{ $A \in s(I)^\circ$, Borel}
\end{equation*}
The original problem of embedding $\mu_X$ in $X$ is equivalent to
the problem of embedding $\mu_Y$ in $Y$ before $Y$ exits the
interval $s(I)^\circ$. Since $Y$ is a time change of a Brownian
motion, we need only consider the problem of embedding $\mu_Y$ in
a Brownian motion before exit from $s(I)^\circ$. If $T$ is an
embedding of $\mu$ in $B$ then $T^{X} \equiv \tau^{-1}(T)$ is
simultaneously an embedding of $\mu_Y$ in $Y$ and $\mu_X$ in $X$.

The first question is when does any embedding exist? If we define
$m =\int_I s(x) \mu_X(dx)$, then the following lemma (see
\citet{PedersenPeskir:01} and
\citet{CoxHobson:02}) gives us necessary and sufficient conditions
for an embedding to exist.

\begin{lemma}
There are three different cases:
\begin{enumerate}

\item $s(I)^\circ = \R$, in which case $X$ is recurrent and we can
embed any distribution $\mu_X$ on $I^\circ$ in $X$.

\item $s(I)^\circ = (-\infty, \alpha)$ (respectively $(-\alpha,
\infty)$) for some $\alpha >0$. Then we may embed $\mu_X$ in $X$
if and only if $m$ exists and $m \ge 0$ (resp.\ $m \le 0$).

\item $s(I)^\circ = ( \alpha , \beta)$, $\alpha < 0 < \beta$. Then
we may embed $\mu_X$ in $X$ if and only if $m =0$.
\end{enumerate}
\end{lemma}

In each case it is clear that:
\begin{equation*}
\mbox{$T^X$ is minimal for $X \iff T^X$ is minimal for $Y \iff T$ is
minimal for $B$,}
\end{equation*}
where $T= \tau(T^X)$.
Further, since $\mu= \mu_Y$ is concentrated on
$s(I)^\circ$, if the stopping time $T$ is minimal then $T$ will
occur before the Brownian motion leaves $s(I)^\circ$ (this is a
consequence of Theorem~\ref{thm:mainB} in case $(ii)$ and
Theorem~\ref{thm:monroe} in case $(iii)$), and then $T^X$ will be
less than the first explosion time of $X$.

It is now possible to apply the results of previous sections to
deduce a series of corollaries about embeddings of $\mu_X$ in $X$.
Suppose that $\mu_X$ can be embedded in $X$ or equivalently that
$\mu_Y$ can be embedded in $B$, before the Brownian motion leaves
$s(I)^\circ$. Suppose further that in the recurrent case where 
$s(I)^\circ = {\mathbb R}$ the law $\mu = \mu_Y$ is integrable. Let 
$T_{max}$ 
and $T^h_{mod}$ be the optimal embeddings of $\mu$ in $B$ as defined in
Sections~\ref{sec:maxminembed} and \ref{sec:maxmod}. (Observe that
from now on we make the dependence of $T^h_{mod}$ on $h$ explicit
in the notation.) Then we can define $T^X_{max}$ and
$T^{X,h}_{mod}$ by
\begin{equation*}
  T^X_{max} = \tau^{-1} \circ T_{max} \hspace{10mm}
  T^{X,h}_{mod} = \tau^{-1} \circ T^h_{mod}.
\end{equation*}

\begin{cor}
$T^X_{max}$ is optimal in the class of minimal embeddings of
$\mu_X$ in $X$ in the sense that it maximises
\begin{equation*}
  \Pr\left(\max_{t \le T} X_t \ge y\right)
\end{equation*}
uniformly in $y \ge 0$.
\end{cor}

\begin{cor}
$T^{X,h}_{mod}$ is optimal in the class of minimal embeddings of
$\mu_X$ in $X$ in the sense that it maximises
\begin{equation*}
  \Pr\left(\max_{t \le T} (h \circ s)(X_t) \ge y\right)
\end{equation*}
uniformly in $y \ge 0$.
\end{cor}

\begin{cor}
$T^{X,|s^{-1}|}_{mod}$ is optimal in the class of minimal
embeddings of $\mu_X$ in $X$ in the sense that it maximises
\begin{equation*}
  \Pr\left(\max_{t \le T} |X_t| \ge y\right)
\end{equation*}
uniformly in $y \ge 0$.
\end{cor}

\bibliography{general}
\bibliographystyle{jmr}
\end{document}